\pdfoutput=1
\documentclass[a4paper,11pt]{article}

\usepackage{fullpage}
\usepackage{amsmath}
\usepackage{booktabs}
\usepackage{multirow}
\usepackage[pdftex]{graphicx}
\usepackage{authblk}
\usepackage{natbib}
\usepackage{palatino}
\usepackage{textcomp}

\usepackage{latexsym, amssymb} 
\usepackage{bbm, bm, nicefrac} 

\usepackage{mathtools} 
\usepackage{eurosym}

\usepackage[squaren]{SIunits}

\usepackage[english]{babel} 
\addto\captionsenglish{}

\usepackage[T1]{fontenc}
\usepackage{ae,aecompl}
\usepackage[utf8]{inputenc} 

\usepackage[dvipsnames, cmyk]{xcolor}
\definecolor{newgreen}{rgb}{0.22, 0.55, 0.11}
\definecolor{newblue}{rgb}{0.55, 0.86, 1}

\usepackage{float, caption, subcaption, colortbl}

\usepackage{hyperref}
\hypersetup{
    colorlinks=true,
    linkcolor=black,
    citecolor=black,
    filecolor=black,
    urlcolor=black,
}

\newcommand{\R}{\mathbbm{R}}
\newcommand{\N}{\mathbbm{N}}

\newcommand{\e}{\varepsilon}
\renewcommand{\d}{x^f}

\newcommand{\ebar}{\overline{\varepsilon}}

\newcommand{\usgwi}{\gram_{\text{CO}_2\text{-eq.}}\per\kilo\watt\hour}
\newcommand{\ugwi}{\kilo\ton_{\text{CO}_2\text{-eq.}}\per\text{a}}
\newcommand{\U}{\mathcal{U}}

\renewcommand{\k}{k}
\newcommand{\K}{\mathcal{K}}

\renewcommand{\cap}{\mathord{\mathit{INVEST}}}

\newcommand{\gwi}{\mathord{\mathit{GWI}}}
\newcommand{\tac}{\mathord{\mathit{TAC}}}
\newcommand{\chp}{\mathord{\mathit{CHP}}}
\newcommand{\pvf}{\mathord{\mathit{PVF}}}

\newcommand{\px}{\chi}
\newcommand{\pp}{\rho}
\renewcommand{\one}{one }
\newcommand{\two}{two }
\newcommand{\three}{three }

\newcommand{\scen}{\xi}
\newcommand{\NAME}{flex-hand\ }
\newcommand{\NAMEBR}{flex-hand}
\newcommand{\NAMEM}{flex\text{-}hand}
\newcommand{\FULLNAME}{flexible here-and-now decision\ }

\title{Flexible here-and-now decisions\\ for two-stage multi-objective
optimization:\\ Method and application to energy system design selection}

\author[1\footnote{dinah.hollermann@rwth-aachen.de}]{Dinah Elena Hollermann}
\author[2\footnote{Corrsponding author: marc.goerigk@uni-siegen.de}]{Marc Goerigk}
\author[1\footnote{doerthe.hoffrogge@rwth-aachen.de}]{D\"orthe Franzisca Hoffrogge}
\author[1\footnote{maike.hennen@rwth-aachen.de}]{Maike Hennen}
\author[1,3\footnote{andre.bardow@ltt.rwth-aachen.de}]{Andr\'e Bardow}

\affil[1]{Institute of Technical Thermodynamics, RWTH Aachen University, 52056 Aachen, Germany}
\affil[2]{Network and Data Science Management, University of Siegen, 57072 Siegen, Germany}
\affil[3]{Institute of Energy and Climate Research, Energy Systems Engineering (IEK-10),
Forschungszentrum J\"ulich GmbH, 52425 J\"ulich, Germany}

\begin{document}

\maketitle

\section*{Abstract}
The synthesis of energy systems is a two-stage optimization problem where design decisions have to be implemented here-and-now (first stage), while for the operation of installed components, we can wait-and-see (second stage). To identify a sustainable design, we need to account for both economical and environmental criteria leading to multi-objective optimization problems. However, multi-objective optimization leads not to one optimal design but to multiple Pareto-efficient design options in general. Thus, the decision maker usually has to decide manually which design should finally be implemented.

In this paper, we propose the \FULLNAME (\NAMEBR) approach for automatic identification of \one single design for multi-objective optimization. The approach minimizes the distance of the Pareto front based on \one fixed design to the Pareto front allowing multiple designs. Uncertainty regarding parameters of future operations can be easily included through a robust extension of the {\NAMEBR} approach.

Results of a real-world case study show that the obtained design is highly flexible to adapt operation to the considered objective functions. Thus, the design provides an energy system with the ability to adapt to a changing focus in decision criteria, e.\,g., due to changing political aims.

\paragraph{Keywords}Multi-objective optimization; automatic solution selection; energy system design; two-stage optimization; robust optimization

\section{Introduction}
\label{sec: intro}
The design of sustainable energy systems needs to balance multiple objectives representing economical, environmental, and social decision criteria. The resulting design problem is therefore best addressed by multi-objective optimization. However, multi-objective optimization yields not an optimal design but a Pareto front with many different designs. Hence, the decision maker is often confronted with the question: How to select \one single design?

In literature, several approaches exist to reduce the number of relevant solutions, so that the decision maker has then to choose \one of fewer options. A method to focus on the relevant solutions a priori by excluding solutions before the optimization is proposed by \cite{Branke2004} and \cite{Rachmawati2009}. They introduce a preference-based evolutionary approach focusing on calculating ``knee'' regions of the Pareto front. \cite{Hennen2017} focus on Pareto-efficient solutions which are near-optimal with respect to an aggregated criterion that represents the overall set of objective functions. An a posteriori approach to reduce the set of relevant Pareto-efficient solutions is to cluster solutions, e.\,g., based on subtractive clustering \citep{Zio2011}, by $k$-means classification \citep{Taboada2007}, or
by a self-organizing map \citep{Li2009}. \cite{Das1999} focus on relevant Pareto-efficient solutions by evaluating subsets of the objective functions. The approach has been further developed by \cite{Antipova2015}. But still, in all approaches which reduce the number of relevant solutions, the decision maker has to select the finally implemented design.

For further reduction of the Pareto-efficient solutions, \cite{Taboada2007} and \cite{Abubaker2014} propose ranking methods which are based on a prioritization of objective functions by the decision maker. For ranking compromising solutions without explicit prioritization, the methods LINMAP (Linear Programming Technique for Multidimensional Analysis of Preference) \citep{Srinivasan1973}, VIKOR \citep{Duckstein1980, Opricovic2004}, and TOPSIS (technique for order preference by similarity to an ideal solution) \citep{Hwang1981,Chen1992} measure the distance from the ideal point and, in TOPSIS, also from the nadir point. The idea is based on compromise programming, where the solution with objective values ``as close as possible'' to the ideal point is chosen, or ``as far away as possible'' from the worst possible point \citep{Zelany1974}. With the aim to determine key players in social networks, \cite{delaFuente2018} employ eleven methods for automatic solution selection within the set of Pareto optimal solutions. These include, e.\,g., highest hypercube \citep{Beume2009}, consensus \citep{Perez2017}, shortest distance to the ideal point \citep{Padhye2011}, and shortest distance to all points. A review on a posteriori decision making has been recently proposed by \cite{Jing2019}. For general multi-objective problems, the introduced ranking methods are suitable to select one solution. However, for synthesis of energy systems, special characteristics of design optimization are beneficial to select one design.

In particular, design optimization of energy systems is a two-stage optimization problem. Two-stage optimization problems consist of two sets of decision variables: The \emph{here-and-now variables} representing the first stage which need to be fixed in the beginning, and the \emph{wait-and-see variables} representing the second stage which can still be adapted later \citep{Ben-Tal2004}. In energy system optimization, the here-and-now variables correspond to the design of energy systems while the operation is determined by the wait-and-see variables on the second stage. The approaches for solution reduction discussed so far do not take any advantage of the two-stage characteristics of energy systems. Thus, these approaches miss the possibility to choose a first-stage solution which provides high adaptability on the second stage. Exploiting adaptability might be particular important. Since future conditions are not known today, the capability to adapt operation to changing circumstances should be targeted \citep{Shang2005}. Political aims might change, and thus, the importance of economical, environmental, and social aims might change. As a result, the sustainable energy system should provide flexible operation to enable an adaptation to a changing focus within the regarded criteria.

Two-stage characteristics are regarded in design selection by \cite{Mattson2003}. In the proposed approach, the design options need to be discrete. For each design option, the corresponding Pareto front is generated. Afterwards, a Pareto filter is applied simultaneously to all generated fronts deleting all dominated solutions. Finally, design options lying inside a pre-defined region of interest are selected. However, with increasing number of components in the given superstructure, the design options along the Pareto front might change more frequently within the region of interest leading to a higher number of suitable design options to select from. A similar approach is proposed by \cite{Carvalho2012}: Based on discrete design options, the corresponding Pareto front is generated. Design options with the ability to undergo large changes in operation enable higher resilience and are thus favored. \cite{Guo2013} introduce a two-stage optimal planning and design method for combined cooling, heating, and power microgrid systems. On the first stage, the system design is optimized using a generic multi-objective optimization approach. On the second stage, the operational costs are minimized. \cite{Wang2019} additionally regard a feedback from the second stage to the first stage to ensure the accuracy of the planning. However, a systematic or even automatic selection of favored design options is not proposed in either approaches, e.\,g., by applying distance measures to assess the generated Pareto fronts.

All discussed approaches do either take advantage of the two-stage characteristic of energy systems \emph{or} propose a ranking of Pareto-efficient solutions. How to select one single design exploiting the two-stage nature has not been proposed so far to the authors' best knowledge.

In this paper, we propose the \emph{\FULLNAME (\NAMEBR) approach} to identify \one single design which represents the whole Pareto front best without depending on any additional information of the decision maker. For this purpose, we minimize the distance between the Pareto front of the synthesis problem, i.\,e., the Pareto front with changing design options, and the Pareto front induced by \one fixed design. The fixed design leading to the minimal distance is identified by the proposed approach since this design provides a high flexibility regarding the considered objective functions. Thus, the second stage can be well adapted to a changing focus from one to another objective function.

Since not only future aims are uncertain but also future parameter values such as demands or costs, we extend our proposed approach by considering uncertain input data based on scenarios. Uncertainties of input data have been regarded, e.\,g., by \cite{Quintana2017, Tock2015}, and \cite{Lemos2018}, using the sensitivity against uncertainties to assess Pareto-efficient solutions. \cite{Gabrielli2019} propose an approach in which Pareto-efficient designs are assessed by performance indicators measuring the robustness and the cost optimality. The final selection of \one design depends on the target levels which need to be provided by the decision maker. \cite{Ide2016} provide an overview of approaches for one-stage robust multi-objective problems. Robust multi-objective optimization has been applied to energy systems by \cite{Majewski2017a}. \cite{Sun2018} propose a multi-objective discrete robust optimization algorithm to identify a single solution by converting multiple objective functions into \one unified cost function. Until now, taking uncertainties into account when automatically  identifying\one single design for two-stage problems has been an open research question.

For certain and uncertain two-stage problems, our approach allows automatic selection of \one design regarding multiple decision criteria. While we introduce our approach in the context of energy systems, the methodology is general and can be applied to any two-stage multi-objective optimization problem.

The remaining article is structured as follows: In Section \ref{sec: method}, we introduce the problem class as well as the \NAME approach and the extension for uncertain input values. In the following Section \ref{sec: case study}, a real-world case study of an industrial park is introduced and the results are evaluated. We give a summary and conclusions in Section \ref{sec: conclusions}.

\section{The \NAME approach for design selection in multi-objective optimization}
\label{sec: method}

The selection of energy system designs considering multiple criteria is complex, since the Pareto front can contain a diverse set of design options. To help the decision maker implementing \one fixed design which allows flexible operation, we propose the \NAME approach. The approach automatically identifies the best possible design looking even beyond the set of Pareto solutions.
Before presenting the \NAME approach in Section \ref{sec: approach} and the robust extension in Section \ref{sec: rob approach}, we first introduce some basic concepts and notation in Section \ref{sec: basics}.

\subsection{Two-stage multi-objective optimization}
\label{sec: basics}

Two-stage optimization depends on \two stages of decision making. Thus, there are \two sets of variables: $\mathcal{X}_f$, the set of feasible solutions for first-stage variables $x^f$, and $\mathcal{X}_s(x^f)$, the set of feasible solutions for the second-stage variables $x^s$ which depends on the chosen first-stage solution. First-stage variables $x^f$ are also called \emph{here-and-now variables} since once these variables are fixed, they cannot be adapted later. Second-stage variables $x^s$ are also called \emph{wait-and-see variables} since they can still be adapted later. As an example, first-stage variables $x^f$ may model an investment in heating equipment (a design decision), while second-stage variables $x^s$ models the way the equipment is run (an operational decision).

In this paper, we use multiple objectives to design sustainable energy systems. As objectives might be conflicting, we are interested in a set of trade-off solutions. A solution is called \emph{Pareto efficient} if there is no other solution that is at least as good in each objective and strictly better in at least one objective \citep{Ehrgott2005a}. The two-stage multi-objective problem with $K$ objectives is given by:
\begin{align*}
\min\ &\Big(f_1(x^f,x^s),\ldots,f_K(x^f,x^s)\Big) \\
\text{s.t. } & x^f \in \mathcal{X}_f \\
& x^s \in \mathcal{X}_s(x^f)\,.
\end{align*}
We assume the set of Pareto-efficient solutions to be discrete and denote them as $(\px^f_1,\px^s_1)$ $,\ldots,$ $(\px^f_N,\px^s_N)$. The set of Pareto-efficient solutions in the objective space, called \emph{Pareto front}, is denoted as $\mathcal{P}^* = ( \pp^1, \ldots, \pp^N )$ with $\pp^i = (\pp^i_1,\ldots,\pp^i_K) = \Big(f_1(\px^f_i,\px^s_i),\ldots,f_K(\px^f_i,\px^s_i)\Big)\in \mathbb{R}^K$. In this paper, we assume that the Pareto front is discrete; for continuous Pareto fronts, our approach is still applicable to a discrete representative set of the efficient solutions.

We call the problem \emph{ideal} when both first-stage variables $x^f$ and second-stage variables $x^s$ can be chosen separately for each efficient point. The corresponding set of Pareto-efficient solutions in the objective space $\mathcal{P}^*$ is called \emph{ideal Pareto front}. The word \emph{ideal} emphasizes that the ideal Pareto front is always better than the Pareto front with fixed first-stage variables.
In energy system optimization, the ideal Pareto front would imply changing design options (e.\,g. heating equipment) along the Pareto front and thus, cannot be reached by energy systems implemented in the real world. In our \NAME approach, we use the ideal Pareto front as benchmark for evaluating Pareto fronts with fixed design options, i.\,e., with fixed first-stage variables.

\subsection{The \FULLNAME approach}
\label{sec: approach}
The idea of the \NAME approach is to find \one fixed design which represents the ideal Pareto front of the synthesis problem best. In optimization of energy systems, the \emph{operation} of an installed system can be adapted but changing the installed system \emph{design} is not possible in the short term. The design with the highest flexibility in operation regarding the objectives is chosen by the \NAME approach. To determine the highest flexibility, the \NAME approach minimizes the distance between the ideal Pareto front and the Pareto front based on \one fixed design but flexible operation.

The \NAME approach is not limited to energy system optimization but can be applied to any two-stage multi-objective problem where the first-stage variables $x^f$ need to be fixed right now and the second-stage variables $x^s$ can be determined later.

For \one fixed design, i.\,e., for a given first-stage solution $x^f$, we calculate Pareto-efficient solutions:
\begin{align*}
\min\ &\Big(f_1(x^f,x^s),\ldots,f_K(x^f,x^s)\Big)\\
\text{s.t. } & x^s \in \mathcal{X}_s(x^f)\,. 
\end{align*}
We obtain a Pareto front $\mathcal{P}(x^f) = \Big(p^1(x^f),\ldots,p^{N(x^f)}(x^f)\Big)$ depending on the first-stage solution $x^f$ with $N(x^f)$ points which we call \emph{fixed first-stage Pareto front}.

Now the question arises: How to choose first-stage variables $x^f$ such that we get the ``best'' design? To determine the quality of a first-stage solution $x^f$, we compare the Pareto front $\mathcal{P}(x^f)$ with fixed first-stage to the ideal Pareto front $\mathcal{P}^*$. For the comparison of \two sets of multi-objective solutions, a variety of distance measures have been developed (for a review see \cite{Zitzler2003}). Here, we choose a comparison metric based on an additive \emph{binary $\e$-indicator}. As discussed by \cite{Zitzler2003}, there is no single best way to compare Pareto fronts, but the $\e$-indicator is proposed as a good overall method. Employing other metrics would be possible in our setting.

For \two sets $\mathcal{P}^1=(p^1,\ldots, p^S) $ and $\mathcal{P}^2 = (q^1,\ldots,q^T)$ in the $K$-dimensional objective space, the binary $\e$-indicator is obtained by
\[ I(\mathcal{P}^1,\mathcal{P}^2) = \min\Big\{ \e\ :\ \forall l\in[T]\ \exists j\in[S] \text{ s.t. } p^j_i -q^l_i \le \e\ \forall i\in[K] \Big\} \]
where we use the notation $[Z] = \{1,\ldots,Z\}$ for sets with any integer $Z\in \N$. For our approach, the measure $I\Big(\mathcal{P}(x^f),\mathcal{P}^*\Big)$ indicates the distance between a fixed first-stage Pareto front and the ideal Pareto front.

\begin{figure}[H]
\centering
\includegraphics[width=9cm]{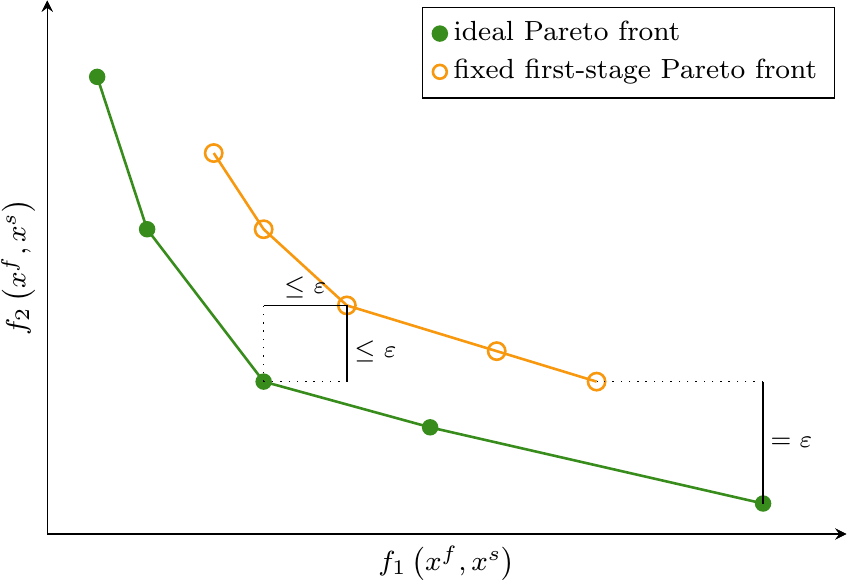}
\caption{Comparison of ideal Pareto front to an arbitrary fixed first-stage Pareto front, to assess the the quality of the fixed first-stage solutions $x^f$; dark green dots: ideal Pareto front used as benchmark; orange circles: fixed first-stage Pareto front with distance $\e$ to the ideal Pareto front; lines are included to guide reader's eyes}
\label{fig: example_nom}
\end{figure}

The comparison metric can be interpreted as follows: Recall that each point of the ideal Pareto front $\mathcal{P}^*$ may involve changing first-stage solutions $x^f$ along the front. In contrast, the fixed first-stage Pareto front $\mathcal{P}(x^f)$ is based on \one single first-stage solution $x^f$ where only the second-stage decision can be adapted. Figure \ref{fig: example_nom} shows the comparison of the ideal Pareto front to an arbitrary fixed first-stage Pareto front.

For each point on the ideal Pareto front, a point on the fixed first-stage Pareto front can be determined such that the difference in each objective function is smaller than a value of $\e$. By minimizing $\e$, the distance between the Pareto fronts is minimized. In this way, for energy systems, we minimize the distance between the ideal Pareto front allowing changing designs as well as operation and the Pareto front based on one fixed design where only the operation can be adapted.

Since we consider the difference regarding each objective function separately, we normalize the objective functions based on the value range in $\mathcal{P}^*$ to circumvent misleading effects by different scales of objective values. To this end, we use the normalized objectives $\overline{f}_i$ with
\[ \overline{f}_i(x^f,x^s)=\frac{ f_i(x^f,x^s)-\min_{j\in [N]}{ f_i(\px^f_j,\px^s_j)}}{\max_{j\in [N]}{ f_i(\px^f_{j},\px^s_{j})}-\min_{j\in [N]}{ f_i(\px^f_{j},\px^s_{j})}}\,.\]
In the following, we write $\overline{\mathcal{P}}(x^f)$, $\overline{\mathcal{P}}^*$ to denote normalized Pareto fronts.

Searching for an optimal first-stage variable $x^f$, i.\,e., an optimal design, we now minimize the distance between the Pareto fronts:
\[ \min\left\{ I\left(\overline{\mathcal{P}}(x^f), \overline{\mathcal{P}}^*\right) : x^f \in \mathcal{X}_f \right\}\,. \]
The problem formulation can also be written as:
\begin{align*}
\min &\ \ebar \\
\text{s.t. } & \overline{f}_i(x^f,x^s_j) - \overline{f}_i(\px^f_j,\px^s_j) \le \ebar & \forall i\in[K],j\in[N] \\
& x^f \in \mathcal{X}_f \\
& x^s_j \in \mathcal{X}_s(x^f) & \forall j\in[N]\,.
\end{align*}
Here, for each point on the ideal Pareto front $\mathcal{P}^*$, we consider \one point on the fixed first-stage Pareto front $\overline{\mathcal{P}}(x^f)$. Thus, for the purpose of finding an optimal $x^f$, we assume without loss of generality that the number of points on the fixed first-stage Pareto front is identical to the number of points on the ideal Pareto front and thus $N(x^f) = N$ holds.

The \NAME approach yields an optimal first-stage solution $(x^f)^*$ which represents the ideal Pareto front best regarding the chosen measure. We call the optimal first-stage solution $(x^f)^*$ the \emph{\NAME solution}. In general, the \NAME solution chosen by our approach is not necessarily part of the solutions of the calculated ideal Pareto front $\mathcal{P}^*$. Thus, approaches based on sorting or solution-reduction would not identify the \NAME solution in general.

Having found a \NAME solution, we calculate the corresponding Pareto-efficient second-stage solution $x^s$ in a separate post-processing optimization step. The resulting Pareto front is called the \emph{\NAME Pareto front}. The number of points on the \NAME Pareto front might differ from the original number of points $N$. For energy systems, this post-processing optimization corresponds to an operational multi-objective optimization based on a fixed design.

There is also an alternative interpretation of our proposed \NAME approach: To solve multi-objective optimization problems, weighted sums as $\sum_{i\in[K]} \lambda_i \overline{f}_i(x^f,x^s)$ could be considered (see \cite{Ehrgott2005a}). Since we do not know the correct preference weighting, the weighting factors $\lambda_i$ are uncertain. Here, each point on the ideal Pareto front represents the optimal solutions if we knew the preference weighting in advance. The aim is now to find a first-stage solution $x^f$ which yields a high solution quality for a wide range of weights. The $\ebar$ value of our \NAME approach can thus be considered as regret (see \cite{aissi2009min} for a survey of regret optimization). However, employing weighted sums, only solutions on the convex hull of the Pareto front can be found. In contrast, our approach also considers points not on the convex hull of Pareto solutions.

\subsection{The robust \NAME approach}
\label{sec: rob approach}

In optimization of energy systems, decisions are based on input parameters which are inherently uncertain. Thus, we extend the proposed \NAME approach for problems comprising uncertain parameters in the objective functions and constraints, and introduce the \emph{robust \NAME approach}. As the \NAME approach, the robust \NAME approach can also be applied to any other two-stage multi-objective problem where the first-stage needs to be determined in advance.

The robust \NAME approach automatically selects first-stage solutions taking uncertainties into account. For this purpose, we assume multiple scenarios which are contained in the discrete uncertainty set $\U$. In each scenario, we compare the ideal Pareto front for the current scenario to a fixed first-stage Pareto front which is based on \one fixed set of first-stage variables for \emph{all} scenarios simultaneously. We then minimize the distance between the Pareto fronts in the worst case and thereby find the robust optimal first-stage solution.

For this purpose, we calculate the ideal Pareto front $\mathcal{P}^*(\xi)$ in each scenario $\xi\in\U$ separately. Here, the number of elements in $\mathcal{P}^*(\xi)$ is denoted by $N(\xi)$. The objectives are parametrized also through scenarios $\xi\in\U$, i.\,e., we use $f_i\Big(x^f,x^s(\xi),\xi\Big)$. Again, we normalize objectives which we denote by $\overline{f}_i\Big(x^f,x^s(\xi),\xi\Big)$, for each scenario $\xi$. To minimize the worst-case distance between the ideal Pareto front and the fixed first-stage Pareto front, we minimize the maximum value of the $\e$-indicator over all $\xi\in\U$:

\begin{align*}
\min &\ \ebar \\
\text{s.t. } & \overline{f}_i\Big(x^f,x^s_j(\xi),\xi\Big) - \overline{f}_i\Big(\px^f_j(\xi),\px^s_j(\xi),\xi\Big) \le \ebar & \forall \xi\in\U, i\in[K],j\in[N(\xi)] \\
& x^f \in \mathcal{X}_f \\
& x^s_j(\xi) \in \mathcal{X}_s(x^f,\xi) & \forall \xi\in\U, j\in[N(\xi)]\,.
\end{align*}
Here, $\mathcal{X}_s(x^f,\xi)$ is the set of feasible second-stage variables given first-stage variables $x^f$ and scenario $\xi$. The optimal first stage solution $(x^f)^*$ is the identified robust \NAME solution.

\begin{figure}[H]
\centering
\includegraphics[width=14cm]{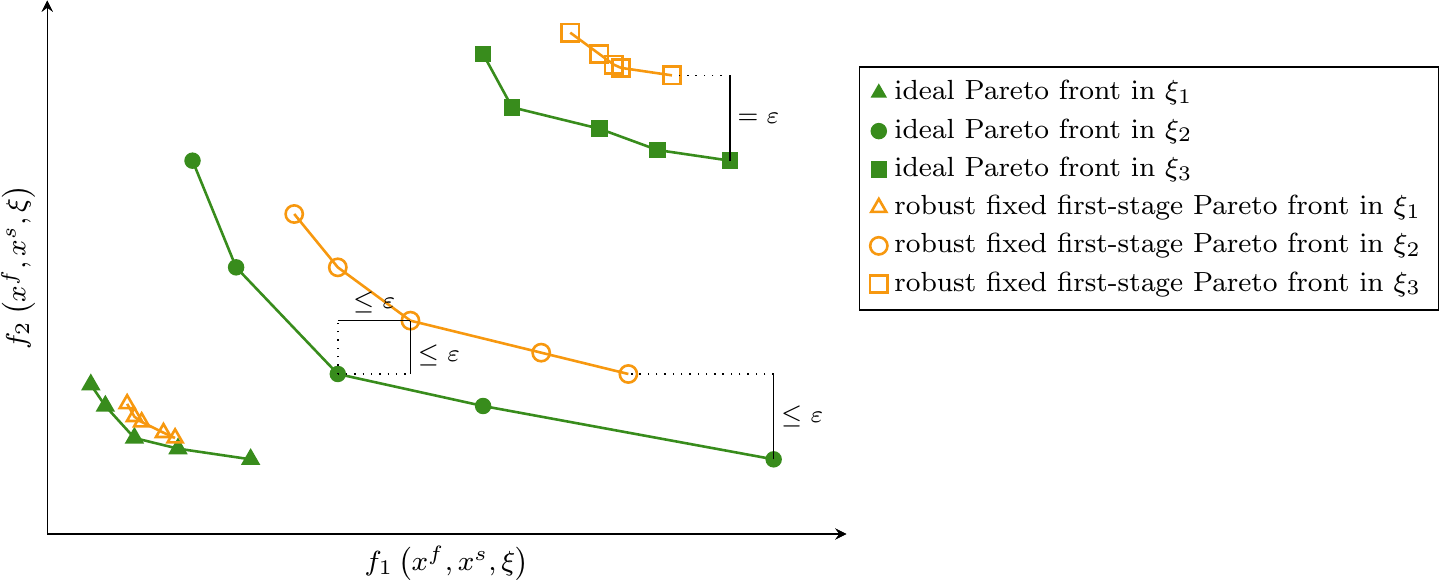}
\caption{Idea of the robust \NAME approach: For each scenario separately, the ideal Pareto front (dark green filled marks) and a robust fixed first-stage Pareto front (orange unfilled marks) are compared; triangles, circles, and squares represent scenario $\scen_1$, $\scen_2$, and $\scen_3$, respectively; the distance $\e$ is calculated regarding all scenarios; here, Pareto fronts are presented without normalization; lines are included to guide reader's eyes}
\label{fig: example_unc}
\end{figure}

An example is given in Figure \ref{fig: example_unc}. The ideal Pareto front is calculated for each of the three scenarios $(\xi_1, \xi_2,\xi_3)$ separately. The robust fixed first-stage Pareto fronts are based on \one fixed set of first-stage variables for all scenarios; however, the corresponding robust fixed first-stage Pareto fronts are calculated for each scenario separately by adapting the second-stage variables $x^s$.

\section{Case study}
\label{sec: case study}

In this section, we apply the \NAME approach to design the energy system of a real-world industrial park. The case study is introduced in Section \ref{sec: real world}. The results of the \NAME approach are presented and discussed in detail in Section \ref{sec: results DS} and in Section \ref{sec: results rob DS} for the robust \NAME approach.

To compute Pareto fronts, we use the adaptive normal boundary intersection method \citep{Das1998}. For calculation, we employ 4 threads of a computer with 3.24 GHz and
64 GB RAM. The problem is formulated in GAMS 24.7.3 \citep{McCarl2016} and solved by the solver CPLEX 12.6.3.0 \citep{IBMILOG2015} to machine accuracy.

\subsection{The real-world example}
\label{sec: real world}

\begin{figure}[H]
\centering
	\includegraphics[width=9cm]{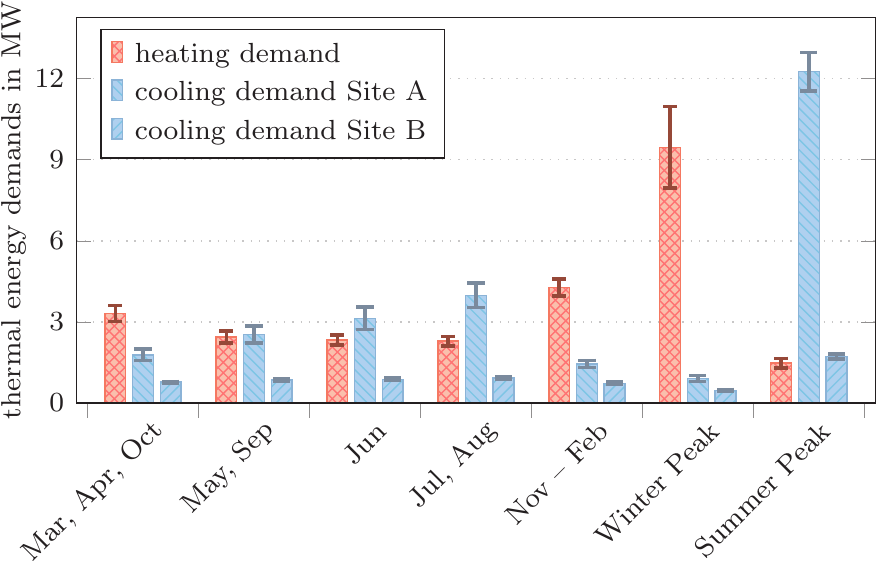}
\caption{Thermal demands of the industrial site and their uncertainties represented by error bars; adopted from \cite{Majewski2017a}}
\label{fig: demands}
\end{figure}

The real-world example is based on our previous work \citep{Voll2013c} on the optimization of a distributed energy supply system. We consider an industrial site with \one power grid, \one heating grid and \two separated cooling grids (Site A and Site B). The thermal demands and their uncertainties are given in Fig.\,\ref{fig: demands}.
The design of the energy system corresponds to sizing and installing any number of components from the following types of energy conversion components: boilers $B$, combined heat and power engines $\chp$, absorption chillers $AC$, and compression chillers $CC$. Natural gas can be used at costs of $p^{gas}=6\,\text{ct}\per\kilo\watt\hour$ with $\pm 40\,\%$ of uncertainty. Furthermore, we assume a connection to the electricity grid. Electricity can be purchased for $p^{el,buy}=16\,\text{ct}\per\kilo\watt\hour$ and sold for $p^{el,sell}=10\,\text{ct}\per\kilo\watt\hour$. For purchasing and selling, an uncertainty of $\pm 46\,\%$ is considered. All possibly uncertain input parameters depend on the scenario $\xi$ and are additionally marked using a tilde. The values for uncertainty are deduced from \cite{Majewski2017a}. To design a sustainable energy system, we employ an economical and an environmental objective function: the total annualized costs $\tac$ and the global warming impact $\gwi$. In principle, the method could also consider social criteria \citep{Mota2015}.

The total annualized costs $\tac$ are defined by:
\begin{align}
&\tac\left(\dot{U},\dot{U}^{el,buy},\dot{V}^{el,sell},\mathord{\mathit{INVEST}}_k; \xi\right)&\\
&\qquad =\sum_{t\in [T]} \Big[\Delta \tau_t \Big( \widetilde{p}^{gas}(\xi) \cdot \sum_{k\in B \cup \chp}\dot{U}_{kt} (\xi) + \widetilde{p}^{el,buy}(\xi) \cdot \dot{U}^{el,buy}_t(\xi) \nonumber\\
& \hspace{6.8cm}- \widetilde{p}^{el,sell}(\xi) \cdot \dot{V}^{el,sell}_t(\xi)\Big)\Big]\nonumber\\
& \hspace{1.9cm}+ \sum_{{\k} \in \K} \bigg(\frac{1}{\pvf}+p_{\k}^{m}\bigg) \cdot \mathord{\mathit{INVEST}}_k\nonumber
	\end{align}
%
Here, $\k$ represents a component in the set of all components $\K=B \cup \chp \cup AC \cup CC$ which might be installed. For each time step $t\in [T]$, $\Delta \tau_t$ represents its length. The corresponding input energy flows of natural gas for boilers $B$ and combined heat and power engines $\chp$ are denoted by $\dot{U}_{kt}$. The input and the output energy flow of electricity are declared by $\dot{U}^{el,buy}_t$ and $\dot{V}^{el,sell}_t$, respectively. For each component ${\k}$, $p_{\k}^{m}$ represents the  annual maintenance costs as share of the investment costs $\cap_{\k}$. For annualizing the investment costs $\cap_{\k}$, we use the present value factor \citep{Broverman2010}
\[\pvf= \frac{(i+1)^h-1}{(i+1)^h\cdot i}\]
with an interest rate $i=8\,\%$ and a time horizon $h=4\,\text{a}$.

The global warming impact is given by:
\begin{align}
&\mathord{\mathit{GWI}}\left(\dot{U},\dot{U}^{el,buy},\dot{V}^{el,sell}; \xi\right)\label{eq: gwi}\\
&\qquad = \sum_{t \in [T]} \Delta \tau_t \left[\sum_{k\in B \cup \chp}\dot{U}_{kt} (\xi)\cdot \mathord{\mathit{GWI}}^{gas}
	+ \left(\dot{U}_t^{el,buy}(\xi) - \dot{V}_t^{el,sell}(\xi)\right)\cdot\widetilde{\mathord{\mathit{GWI}}}^{el}(\xi)\right]\,.\nonumber
	\end{align}
We employ $\mathord{\mathit{GWI}}^{gas}=244\,\usgwi$ for the specific global warming impact of gas which is not subject to remarkable variation. For the specific global warming impact $\mathord{\mathit{GWI}}^{el}$ of electricity purchased from the grid, we employ a value of $561\,\usgwi$. Since the future electricity mix might change significantly, we assume the specific global warming impact $\mathord{\mathit{GWI}}^{el}$ to be uncertain lying within $430\,\usgwi$ and $610\,\usgwi$. When selling electricity to the grid, a credit for global warming impact is given, following the idea of the avoided burden \citep{Baumann2004}. Here, the global warming impact $\gwi$ depends only implicitly on the first-stage variables $\d$ due to the constraints. A direct influence would be given if the global warming impact induced by the manufacturing of the components was taken into account. However, since the global warming impact of the operation has usually a significantly higher impact \citep{Guillen-Gosalbez2011}, we neglect this dependency. The complete \NAME optimization model is presented in Appendix~\ref{sec: appendix}.

For the design optimization, we assume a ``green field" without existing energy components. However, the \NAME approach could also be applied to retrofit an energy system.

\subsection{The \NAME design}
\label{sec: results DS}

\begin{figure}[H]
\centering
	\includegraphics[width=9cm]{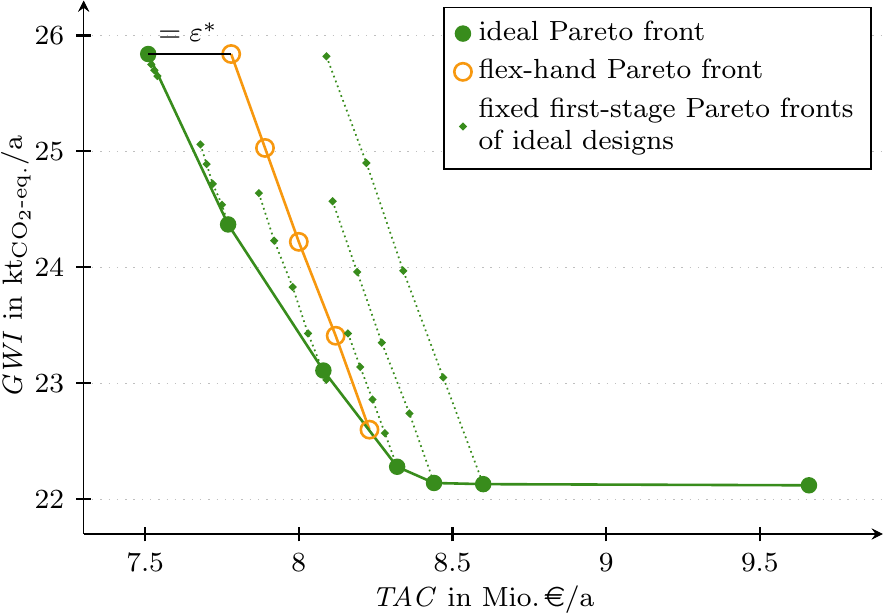}
\caption[Comparison of ideal Pareto front and the \NAME Pareto front]{Comparison of ideal Pareto front (dark green dots) and the \NAME Pareto front (orange circles) with minimal distance $\e^*$ to the ideal Pareto front; small dark green dots: fixed first-stage Pareto front of ideal designs, i.\,e., Pareto-efficient operation for each design of the ideal Pareto front; here, Pareto fronts are presented without normalization; lines are included to guide reader's eyes}
\label{fig: cs pf}
\end{figure}

We now employ the \NAME approach to design the sustainable energy system in order to obtain the best solution for the first-stage variables $x^{f}$ which we call the \emph{\NAME design}. For this purpose, we first calculate the ideal Pareto front as a benchmark. The ideal Pareto front is obtained by allowing a different design for each point on the front. The largest optimization problem for calculating a point on the ideal Pareto front consists of $1950$ equations, $576$ variables, and $310$ binary variables after presolve. In total, the whole ideal Pareto front is calculated in $317\,\second$. The \NAME problem has $5099$ equations, $2023$ variables, and $751$ binary variables after presolve. Here, computing the \NAME Pareto front takes $152\,\second$. Both the ideal Pareto front and the \NAME Pareto front are shown in Fig.\,\ref{fig: cs pf}.

The \NAME Pareto front of the selected design is ``streched out''. Thus, the \NAME design allows for flexible operation providing a high ability to adapt to changing future objectives. The \NAME design can be operated such that the total annualized costs $\tac$ are very low at $7.8\,\text{Mio.\euro}\per\text{a}$ or the global warming impact $\gwi$ is very low with $22.6\,\ugwi$.

In this case study, the minimal distance between the ideal and the \NAME Pareto front is limited, e.\,g., by the anchor points with minimal total annualized costs limits (Fig.\,\ref{fig: cs pf}). The corresponding scaled value for the minimized distance is $\overline{\varepsilon}^*=0.128$. For unscaled values, the minimal total annualized costs for the ideal design are $\tac^{ideal}=7.51\,\text{Mio.\euro}\per\text{a}$ and for the \NAME design $\tac^{\NAMEM}=7.78\,\text{Mio.\euro}\per\text{a}$, respectively. Hence, the maximal deviation for total annualized costs is $0.27\,\text{Mio.\euro}\per\text{a}$ which corresponds to a maximal loss of only $3.6\,\%$ compared to the ideal design with minimized total annualized costs. Regarding the minimal global warming impact, the maximal deviation is only $2.17\,\%$. Thus, the flexibility of the \NAME design is very high regarding both objective functions.

\begin{figure}[H]
\centering
	\includegraphics[width=14cm]{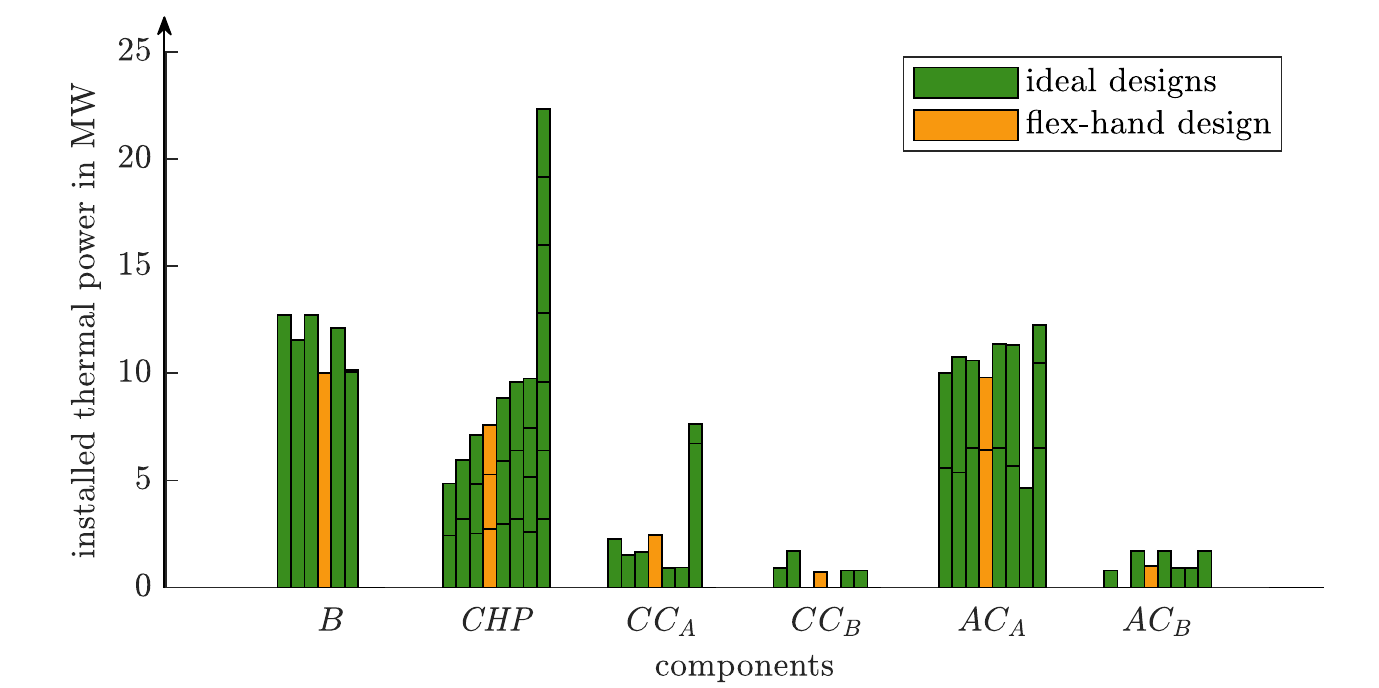}
\caption{In dark green: designs of ideal Pareto front; in orange: \NAME design of the \NAME Pareto front; from left to right, the designs are ordered by decreasing minimal global warming impact; $B$ boiler, $\chp$ combined heat and power engine, $CC_A$ and $CC_B$ compression chillers, and $AC_A$ and $AC_B$ absorption chillers installed on Site A and Site B, respectively}
\label{fig: cs ds}
\end{figure}

When having a closer look at the identified design, we cannot identify a single reason for its higher ability to adapt operation (Fig.\,\ref{fig: cs ds}).
 In general, solutions with lower global warming impact prefer installing higher capacity of combined heat and power engines, since the specific global warming impact of the electricity mix of the grid is higher than the impact of the combined heat and power engines in combination with absorption chillers. The \NAME design does not show remarkable differences compared to the other ideal designs but provides an excellent compromise. Without the proposed approach, this highly adaptable design would most likely not have been identified by the decision maker.

\subsection{The robust \NAME design}
\label{sec: results rob DS}

We now apply the robust \NAME approach to the proposed case study taking uncertainties into account. The uncertainties are introduced in Section \ref{sec: real world}. Here, we consider three scenarios $\scen_1$, $\scen_2$, and $\scen_3$. Scenario $\scen_2$ corresponds to values of the problem without uncertainties discussed in Section \ref{sec: results DS}. In scenario $\scen_1$, we assume all uncertain values to take their smallest values within the uncertainty range and in scenario $\scen_3$ their largest values, respectively. However, any other scenario could be chosen.

		\begin{figure}[H]
\centering
	\includegraphics[width=9cm]{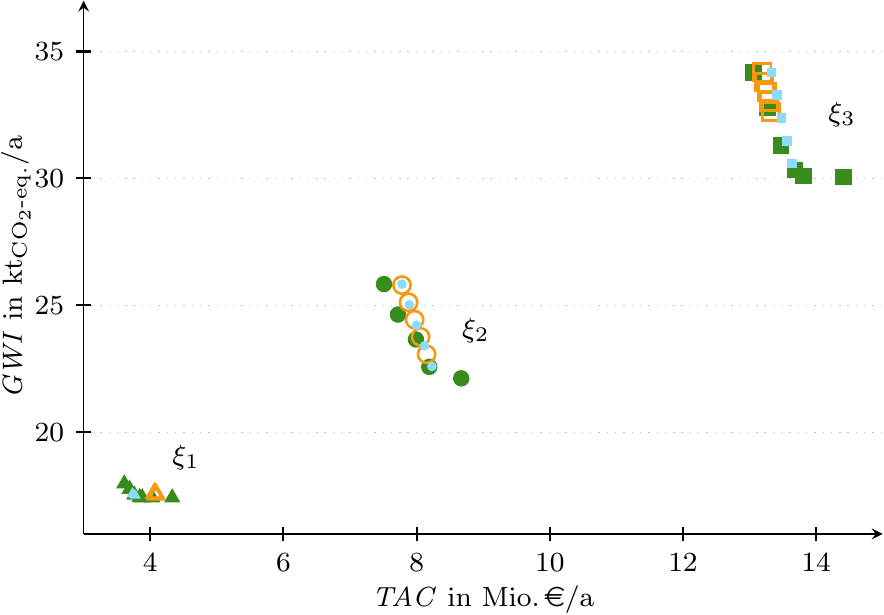}
\caption{Triangles, circles, and squares represent scenario $\scen_1$, $\scen_2$, and $\scen_3$, respectively; dark green filled marks: ideal Pareto front generated for each scenario separately, orange unfilled marks: robust flex-hand Pareto front in each scenario; small light blue marks: \NAME Pareto front for separately considered scenarios; here, Pareto fronts are presented without normalization}
\label{fig: cs rob pf}
\end{figure}

\begin{figure}[H]
\centering
	\includegraphics[width=9cm]{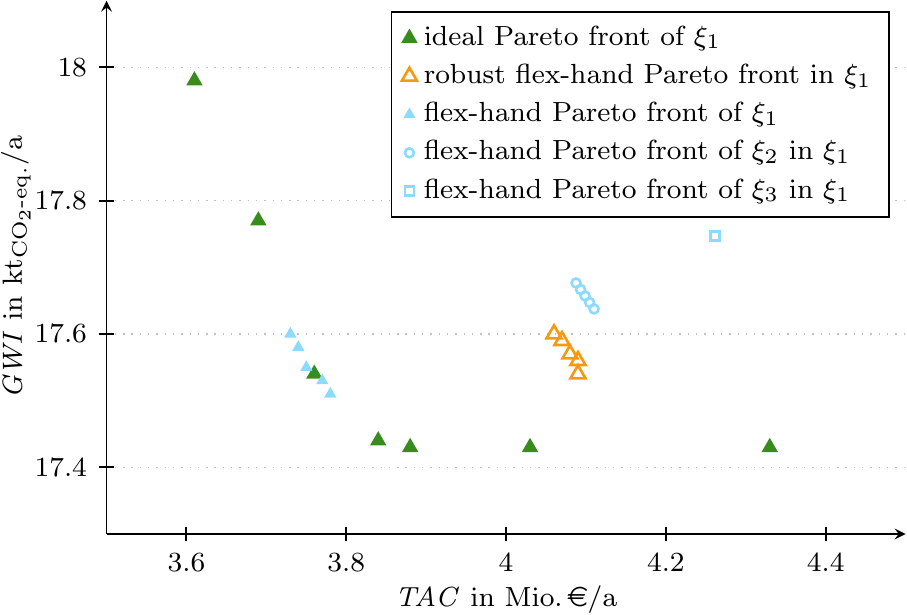}
\caption{All Pareto fronts for scenario $\scen_1$; dark green filled triangles (\textcolor{newgreen}{$\boldsymbol{\blacktriangle}$}): ideal Pareto front of scenario $\scen_1$; orange unfilled triangles (\textcolor{YellowOrange}{$\footnotesize\boldsymbol{\triangle}$}): robust \NAME Pareto front in scenario $\scen_1$; small light blue triangles (\textcolor{newblue}{$\footnotesize\boldsymbol{\blacktriangle}$}): \NAME Pareto front of scenario $\scen_1$; small light blue unfilled circles and squares (\textcolor{newblue}{$\boldsymbol{\circ}$} and \textcolor{newblue}{$\scriptsize\boldsymbol{\square}$}): \NAME Pareto front in scenario $\scen_1$ based on \NAME design computed for scenario $\scen_2$ and $\scen_3$, respectively; here, Pareto fronts are presented without normalization}
\label{fig: ds lb}
\end{figure}

To evaluate the robust \NAME design, we compare the robust \NAME Pareto fronts in all \three scenarios to the \NAME Pareto fronts generated for each scenario separately. Fig.\,\ref{fig: cs rob pf} shows that the \NAME Pareto fronts generated for each scenario separately do not coincide with the robust \NAME Pareto fronts. In scenario $\scen_3$, the robust \NAME design leads to smaller total annualized costs than the \NAME design computed for scenario $\scen_3$ but to a higher global warming impact. In total, the robust \NAME Pareto front is less ``streched out'' than for the nominal case (Section~\ref{sec: results DS}) leading to an optimal distance $\ebar^*$ of $0.625$. The reduced adaptability to the ideal Pareto fronts is due to the fact that the robust \NAME Pareto fronts need to approximate \three ideal Pareto fronts simultaneously, instead of just \one Pareto front. Thus, a good performance of a \NAME design in \one scenario might lead to a poor performance in another scenario if uncertainties are not regarded during design (Fig.\,\ref{fig: ds lb}). In contrast, the robust \NAME design is a compromise solution performing well in all \three scenarios simultaneously.

In Fig.\,\ref{fig: ds lb}, we take a closer look on the computed Pareto fronts in scenario $\scen_1$. Here, the robust \NAME design (\textcolor{YellowOrange}{$\footnotesize\boldsymbol{\triangle}$}) clearly performs better than the \NAME design identified for scenario $\scen_2$ (\textcolor{newblue}{$\boldsymbol{\circ}$}) and the \NAME design identified for scenario $\scen_3$ (\textcolor{newblue}{$\scriptsize\boldsymbol{\square}$}). In scenario $\scen_1$, only the \NAME Pareto front of scenario $\scen_1$ (\textcolor{newblue}{$\footnotesize\boldsymbol{\blacktriangle}$}) approximates the ideal Pareto front (\textcolor{newgreen}{$\boldsymbol{\blacktriangle}$}) better than the robust \NAME Pareto front (\textcolor{YellowOrange}{$\footnotesize\boldsymbol{\triangle}$}). However, the \NAME design of scenario $\scen_1$ is infeasible for scenario $\scen_2$ and $\scen_3$. In contrast, the robust \NAME design is feasible and performs well for all scenarios. 

\begin{figure}[H]
\centering
	\includegraphics[width=14cm]{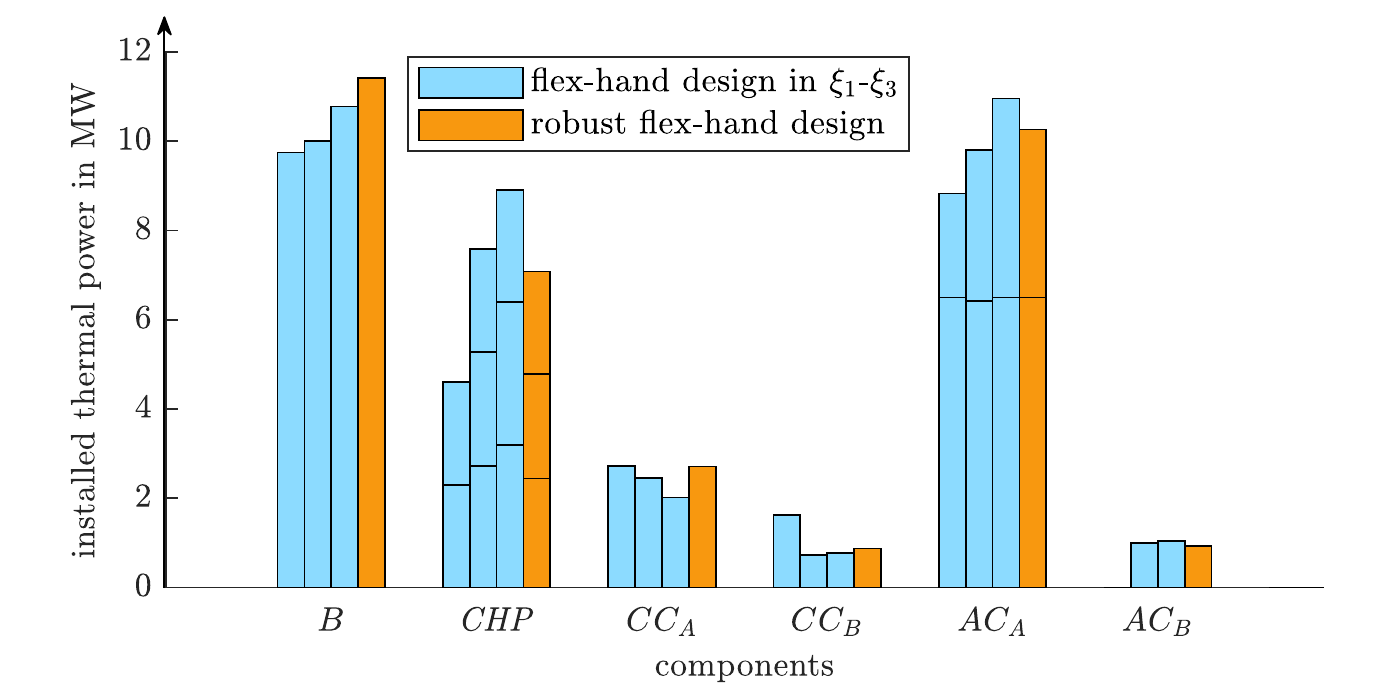}
\caption{In light blue: \NAME designs generated for each scenario separately (from left to right: scenario $\scen_1$, $\scen_2$, and $\scen_3$); in orange: robust \NAME design; $B$ boiler, $\chp$ combined heat and power engine, $CC_A$ and $CC_B$ compression chillers, and $AC_A$ and $AC_B$ absorption chillers installed on Site A and Site B, respectively}
\label{fig: cs rob ds}
\end{figure}

Having a closer look at the design (Fig.\,\ref{fig: cs rob ds}), we observe that the total capacity of the \three \NAME designs increases from scenario $\scen_1$ to $\scen_3$. This is due to the fact that values of uncertain input parameter increase as well. With increasing demands and also increasing specific global warming impact of the electricity grid, larger combined heat and power engines and boilers are installed combined with a higher capacity of absorption chillers and smaller compression chillers. The robust \NAME design does not differ remarkably from the \three \NAME designs. Thus, the robust \NAME approach is necessary to identify the excellent compromise given by the robust \NAME design.

\section{Conclusions}
\label{sec: conclusions}

The sustainable optimization of energy systems is inherently a two-stage optimization problem with multiple decision criteria. Applying multi-objective optimization usually generates different designs for each point on the Pareto front. We propose the \NAME approach for identifying a single design which performs well regarding all decision criteria. The idea of the \NAME approach is to approximate the Pareto front with changing design options (\emph{ideal Pareto front}) by a Pareto front with one fixed design for the whole front. The design leading to the minimal distance between both Pareto fronts is the design which is identified by the \NAME approach. The identified design (\emph{\NAME design}) is able to adapt to all regarded criteria well, and thus, provides high flexibility to reach future aims for which focus might change between the considered objective functions.

Our real-world case study demonstrates the resulting high adaptability with respect to the considered criteria. For designing the sustainable energy system, we choose total annualized costs and the global warming impact as economical and environmental criteria, respectively. The calculated Pareto front of the \NAME design is "stretched out" in comparison to the Pareto fronts obtained by operational optimization of designs lying on the ideal Pareto front. The objective function values of \NAME design differ by less than $3.6\,\%$ from the ideal values which highlights the excellent quality of the identified \NAME design. The \NAME design does not show remarkable differences compared to the designs lying on the ideal Pareto front. Thus, without the \NAME approach, the decision maker would possibly not have chosen the identified design. This effect becomes even more pronounced when considering multiple scenarios simultaneously to account for uncertainty, in which case our approach is able to find a robust solution.

To conclude, the \NAME approach takes advantage of the two-stage nature of energy systems to automatically select \one single design which provides a high flexibility to adapt operation to all considered criteria.

\section*{Acknowledgments}
This work was supported by the Helmholtz Association under the Joint Initiative Energy System 2050---A Contribution of the Research Field Energy.

\appendix

\section{Model Formulation}
\label{sec: appendix}

In the following, we provide the problem formulation of the robust \NAME optimization for the DESS considered in our case study (Section \ref{sec: case study}). In the case study, we consider the total annualized costs $\tac$ and the global warming impact $\gwi$ as objective functions. Uncertainties are regarded for tariffs for purchasing gas $\widetilde{p}^{gas}(\xi)$ and electricity $\widetilde{p}^{el,buy}(\xi)$ as well as for selling electricity $\widetilde{p}^{el,sell}(\xi)$. Furthermore, the specific global warming impact of the electricity mix of the grid $\widetilde{\mathord{\mathit{GWI}}}^{el}(\xi)$  is assumed to be uncertain as well. In the constraints, the energy balances are affected by uncertain energy demands $\widetilde{\dot{E}}^{heat}(\xi),\widetilde{\dot{E}}^{cool}(\xi)$, and $\widetilde{\dot{E}}^{el}(\xi)$.
	
\begin{align*}
\min &\ \ebar \\
\text{s.t. } & \overline{\tac} \left(\dot{U},\dot{U}^{el,buy},\dot{V}^{el,sell},\gamma,\dot{V}^{N}; \xi,j\right) - \overline{\tac}^*(\xi,j) \le \ebar & \forall \xi\in\U,j\in[N(\xi)] \\
& \overline{\gwi} \left(\dot{U},\dot{U}^{el,buy},\dot{V}^{el,sell}; \xi,j\right) - \overline{\gwi}^*(\xi,j) \le \ebar & \forall \xi\in\U, j\in[N(\xi)] \\
&\sum_{k \in B \cup \chp} \dot{V}_{kt} (\xi,j) - \sum_{k\in AC}   \dot{U}_{kt} (\xi,j)= \widetilde{\dot{E}}_t^{heat}(\xi) & \forall t \in [T], \xi\in\U, j\in[N(\xi)] \\
& \sum_{k \in AC \cup CC} \dot{V}_{kt} (\xi,j) = \widetilde{\dot{E}}_t^{cool}(\xi)&\forall t \in [T], \xi\in\U, j\in[N(\xi)]\\
& \sum_{k \in \chp} \dot{V}^{el}_{kt} (\xi,j)- \sum_{k \in CC}  \dot{U}_{kt} (\xi,j)  \\
& \qquad +\ \dot{U}^{el,buy}_t (\xi,j) - \dot{V}_t^{el,sell} (\xi,j) = \widetilde{\dot{E}}_t^{el}(\xi) &  \forall t \in [T],\xi\in\U, j\in[N(\xi)]\\
%
&\sum_{h\in [H]} \gamma_{kh} \leq 1 & \forall k \in \K\\
& \gamma_{kh}  \cdot \dot{V}^{N,lb}_{kh} \leq \dot{V}^{N}_{kh}  \leq \gamma_{kh}  \cdot  \dot{V}^{N,lb}_{kh+1}& \forall\ k \in \K, \forall h \in [H]\\
%
& \rho^{min} \cdot \sum_{h\in[H]} \dot{V}^{N}_{kh}  \leq \dot{V}_{kt} (\xi,j)\leq \sum_{h\in[H]} \dot{V}^{N}_{kh}  	&\forall k \in \K, t \in [T], \xi \in \U, j\in[N(\xi)]\\
&\dot{V}_{kt} (\xi,j) = \eta_k \cdot \dot{U}_{kt} (\xi,j) &\forall k \in \K,t \in [T],\xi\in\U, j\in[N(\xi)]\\
&\dot{V}^{el}_{kt} (\xi,j) = \eta^{tot}_k \cdot \dot{U}_{kt} (\xi,j)-\dot{V}_{kt} (\xi,j) & \forall   k \in \chp,t \in [T],\xi\in\U, j\in[N(\xi)]\\
%
& \ebar \in \R_{+}\\
&\dot{U}^{el,buy}(\xi,j),\dot{V}^{el,sell}(\xi,j),\dot{V}^{el}(\xi,j), 
\dot{U}(\xi,j),\dot{V}(\xi,j)\in \R_{+}^{|\K|\times T}& \forall \xi \in \U, j\in[N(\xi)] \\
&\gamma  \in \{0,1\}^{|\K|\times H},\dot{V}^{N} \in \R_{+}^{|\K|\times H}
%
\end{align*}
The total annualized costs $\tac$ and the global warming impact $\gwi$ are defined by
\begin{align}
&\tac\left(\dot{U},\dot{U}^{el,buy},\dot{V}^{el,sell},\gamma,\dot{V}^{N}; \xi,j\right)&\nonumber\\
&\qquad =\sum_{t\in [T]} \Big[\Delta \tau_t \Big( \widetilde{p}^{gas}(\xi) \cdot \sum_{k\in B \cup \chp}\dot{U}_{kt} (\xi,j) + \widetilde{p}^{el,buy}(\xi) \cdot \dot{U}^{el,buy}_t(\xi,j) \nonumber\\
& \hspace{6.8cm}- \widetilde{p}^{el,sell}(\xi) \cdot \dot{V}^{el,sell}_t(\xi,j)\Big)\Big]\nonumber\\
& \hspace{1.9cm}+ \sum_{{\k} \in \K} \left(\frac{1}{\pvf}+p_{\k}^{m}\right) \cdot \underbrace{\sum_{h\in [H]}\left[ \gamma_{kh} \cdot  \kappa_{kh}+ m_{kh} \cdot \left(\dot{V}^{N}_{kh}  -  \gamma_{kh} \dot{V}^{N,lb}_{kh}\right)\right]}_{\eqqcolon \mathord{\mathit{INVEST}}_k}\nonumber\\
&\mathord{\mathit{GWI}}\left(\dot{U},\dot{U}^{el,buy},\dot{V}^{el,sell}; \xi,j\right)\nonumber\\
&\qquad = \sum_{t \in [T]} \Delta \tau_t \left[\sum_{k\in B \cup \chp}\dot{U}_{kt} (\xi,j)\cdot \mathord{\mathit{GWI}}^{gas}
	+ \left(\dot{U}_t^{el,buy}(\xi,j) - \dot{V}_t^{el,sell}(\xi,j)\right)\cdot\widetilde{\mathord{\mathit{GWI}}}^{el}(\xi)\right]\,.\nonumber
	\end{align}
Bars above the total annualized costs $\tac$ and the global warming impact $\gwi$ in the optimization problem denote the normalization of the objective values. Objective values on the normalized ideal Pareto fronts are denoted by $\Big(\overline{\tac}^*(\xi,j), \overline{\gwi}^*(\xi,j)\Big)$ for each point $j\in [N(\xi)]$ and each scenario $\xi$ in the uncertainty set $\U$.

The duration of a time step $t\in [T]$ is given by $\Delta \tau_t$. Maintenance costs are determined by the share $p_{\k}^{m}$ of the investment costs $\mathord{\mathit{INVEST}}_k$. The investment costs $\mathord{\mathit{INVEST}}_k$ are annualized by the present value factor $\pvf$. $\mathord{\mathit{GWI}}^{gas}$ represent the specific global warming impact of purchased gas. Purchased and sold electricity is denoted by $\dot{U}_t^{el,buy}$ and $\dot{V}_t^{el,sell}$, respectively. $\dot{U}_{kt}$ and $\dot{V}_{kt}$ specifies input and output energy flows in time step $t$ of component $k\in \K$. Components include boilers $B$, combined heat and power engines $\chp$, absorption chillers $AC$, and compression chillers $CC$. Input and output energy flows are coupled by the thermal efficiency $\eta_k$. For combined heat and power engines, the total efficiency $\eta^{tot}_k$ is given by the sum of the thermal and the electrical efficiency $\eta^{tot}_k=\eta_k+\eta^{el}_k$. The minimal part-load of a component $k$ is defined by the fraction $\rho^{min}$ of the installed nominal capacity.

\begin{figure}[H]
\centering
\includegraphics{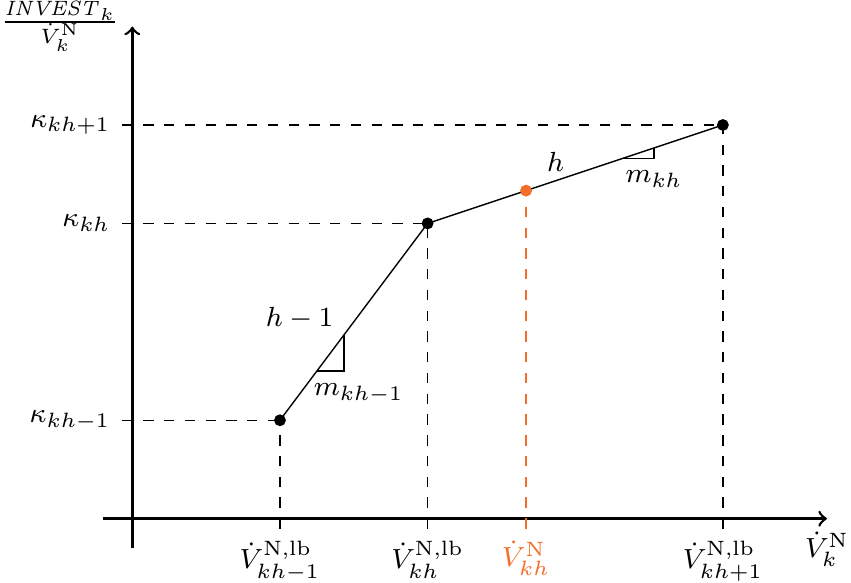}
\caption{Piecewise linearization of the investment costs $\mathord{\mathit{INVEST}}_k$ of a newly installed component $k$ is presented. Here, $h$ is the active line segment; thus, $\gamma_{kh}$ is equal to 1.}
\label{fig: lin invest}
\end{figure}

The investment costs $\mathord{\mathit{INVEST}}_k$ of a newly installed component $k$ are linearized by piecewise linearization with $\sum_{h\in [H]} \left[\gamma_{kh} \cdot  \kappa_{kh}+ m_{kh} \cdot \left(\dot{V}^{N}_{kh}  -  \gamma_{kh} \dot{V}^{N,lb}_{kh}\right)\right]$ (see Fig.\,\ref{fig: lin invest}). $m_{kh}$ is the gradient for each line segment $h\in[H]$ and is defined by
\[m_{kh} \coloneqq \frac{\kappa_{kh+1}-\kappa_{kh}}{\dot{V}^{N, lb}_{kh+1}-\dot{V}^{N,lb}_{kh}}\quad \forall\ k \in \K,  h \in [H]\,.\]
Here, parameters $\dot{V}^{N,lb}_{kh+1}$ and $\dot{V}^{N,lb}_{kh}$ represent the nominal capacities of the lower and upper supporting point of line segment $h$ and parameters $\kappa_{kh}$ and $\kappa_{kh+1}$ the corresponding specific investment costs. Binary variables $\gamma_{kh}$ determine if line segment $h$ is active ($\gamma_{kh}=1$). Since the sum $\sum_{h\in [H]}\gamma_{kh}$ is equal to $1$, only one line segment can be active at the time. Thus, only one value for the nominal capacity $\dot{V}^{N}_{kh}$ of all line segments is unequal to $0$; hence, the nominal capacity $\dot{V}^{N}_{k}$ of an installed component $k$ is given by the sum $\sum_{h\in[H]} \dot{V}^{N}_{kh}$.

{
\bibliographystyle{apalike} 
\bibliography{references}
}

\end{document}